\documentclass[12pt]{article}
\usepackage{amssymb,amsthm,amsmath,amsfonts,latexsym,tikz,hyperref,ytableau}
\usepackage[hmargin=1in,vmargin=1in]{geometry}
\usepackage{shuffle}

\newtheorem{thm}{Theorem}[section]
\newtheorem{prop}[thm]{Proposition}
\newtheorem{cor}[thm]{Corollary}
\newtheorem{lem}[thm]{Lemma}
\newtheorem{conj}[thm]{Conjecture}
\newtheorem{exa}[thm]{Example}
\newtheorem{question}[thm]{Question}

\DeclareMathOperator{\std}{std}
\DeclareMathOperator{\SSYT}{SSYT}
\DeclareMathOperator{\RS}{RS}
\DeclareMathOperator{\MR}{MR}
\DeclareMathOperator{\Pk}{Pk}

\newcommand{\shu}{\shuffle}
\newcommand{\pshu}{\hs{2pt}\ol{\shu}\hs{2pt}}
\newcommand{\ofS}{\ol{\fS}}
\newcommand{\lab}{\ol{\la}}
\newcommand{\cHb}{\ol{\cH}}

\newcommand{\ben}{\begin{enumerate}}
\newcommand{\een}{\end{enumerate}}
\newcommand{\ble}{\begin{lem}}
\newcommand{\ele}{\end{lem}}
\newcommand{\bth}{\begin{thm}}
\renewcommand{\eth}{\end{thm}}
\newcommand{\bpr}{\begin{prop}}
\newcommand{\epr}{\end{prop}}
\newcommand{\bco}{\begin{cor}}
\newcommand{\eco}{\end{cor}}
\newcommand{\bcon}{\begin{conj}}
\newcommand{\econ}{\end{conj}}
\newcommand{\bde}{\begin{defn}}
\newcommand{\ede}{\end{defn}}
\newcommand{\bex}{\begin{exa}}
\newcommand{\eex}{\end{exa}}
\newcommand{\barr}{\begin{array}}
\newcommand{\earr}{\end{array}}
\newcommand{\btab}{\begin{tabular}}
\newcommand{\etab}{\end{tabular}}
\newcommand{\beq}{\begin{equation}}
\newcommand{\eeq}{\end{equation}}
\newcommand{\bea}{\begin{eqnarray*}}
\newcommand{\eea}{\end{eqnarray*}}
\newcommand{\bal}{\begin{align*}}
\newcommand{\bce}{\begin{center}}
\newcommand{\ece}{\end{center}}
\newcommand{\bpi}{\begin{picture}}
\newcommand{\epi}{\end{picture}}
\newcommand{\bpp}{\begin{picture}}
\newcommand{\epp}{\end{picture}}
\newcommand{\bfi}{\begin{figure} \begin{center}}
\newcommand{\efi}{\end{center} \end{figure}}
\newcommand{\bprf}{\begin{proof}}
\newcommand{\eprf}{\end{proof}\medskip}
\newcommand{\capt}{\caption}
\newcommand{\bsl}{\begin{slide}{}}
\newcommand{\esl}{\end{slide}}
\newcommand{\bfr}{\begin{frame}}
\newcommand{\efr}{\end{frame}}

\newcommand{\comp}{\models}

\newcommand{\hqed}{\hfill \qed}

\newcommand{\eqed}[1]{\vs{10pt}$\textcolor{white}{\qed}\hfill{\dil#1}\hfill\qed$}

\newcommand{\ol}{\overline}

\newcommand{\hs}[1]{\hspace{#1}}
\newcommand{\hso}[1]{\hspace{-1pt}}
\newcommand{\vs}[1]{\vspace{#1}}
\newcommand{\qmq}[1]{\quad\mbox{#1}\quad}

\newcommand{\emp}{\emptyset}

\newcommand{\sbs}{\subset}
\newcommand{\sbe}{\subseteq}

\newcommand{\ptn}{\vdash}

\newcommand{\case}[4]{\left\{\barr{ll}#1&\mbox{#2}\\#3&\mbox{#4}\earr\right.}

\def\<{\langle}
\def\>{\rangle}

\newcommand{\ree}[1]{(\ref{#1})}

\newcommand{\ra}{\rightarrow}

\newcommand{\al}{\alpha}
\newcommand{\be}{\beta}

\newcommand{\de}{\delta}

\newcommand{\io}{\iota}
\newcommand{\ka}{\kappa}
\newcommand{\la}{\lambda}

\newcommand{\si}{\sigma}

\newcommand{\Th}{\Theta}

\newcommand{\bx}{{\bf x}}

\newcommand{\bbQ}{{\mathbb Q}}

\newcommand{\cA}{{\cal A}}

\newcommand{\cH}{{\cal H}}

\newcommand{\cP}{{\cal P}}

\newcommand{\cS}{{\cal S}}
\newcommand{\cT}{{\cal T}}

\newcommand{\cZ}{{\cal Z}}

\newcommand{\fS}{{\mathfrak S}}

\newcommand{\Hb}{\ol{H}}

\DeclareMathOperator{\Av}{Av}

\DeclareMathOperator{\des}{des}
\DeclareMathOperator{\Des}{Des}

\DeclareMathOperator{\sh}{sh}

\DeclareMathOperator{\Sym}{{Sym}}
\DeclareMathOperator{\QSym}{{QSym}}
\DeclareMathOperator{\SYT}{SYT}

\newcommand{\dil}{\displaystyle}

\begin{document}
\pagestyle{plain}

\title{Pattern avoidance and quasisymmetric functions
}
\author{Zachary Hamaker\\[-5pt]
\small Department of Mathematics, University of Michigan,\\[-5pt]
\small Ann Arbor, MI 48109-1043, {\tt zachary.hamaker@gmail.com}\\
Brendan Pawlowski\\[-5pt]
\small Department of Mathematics, University of Southern California\\[-5pt]
\small Los Angeles, CA 90089-2532, {\tt bpawlows@usc.edu}\\
Bruce E. Sagan\\[-5pt]
\small Department of Mathematics, Michigan State University,\\[-5pt]
\small East Lansing, MI 48824-1027, {\tt sagan@math.msu.edu}
}

\date{\today\\[10pt]
	\begin{flushleft}
	                                       \\[5pt]
	\end{flushleft}}

\maketitle

\begin{abstract}
Given a set of permutations $\Pi$, let $\fS_n(\Pi)$ denote the set of permutations in the symmetric group $\fS_n$ that avoid every element of $\Pi$ in the sense of pattern avoidance.  Given a subset $S$ of $\{1,\dots,n-1\}$, let $F_S$ be the fundamental quasisymmetric function indexed by $S$.  Our object of study is the generating function $Q_n(\Pi) =\sum F_{\Des\si}$ where the sum is over all $\si\in\fS_n(\Pi)$ and $\Des\si$ is the descent set of $\si$.  We  characterize those $\Pi\sbe\fS_3$ such that $Q_n(\Pi)$ is symmetric or Schur nonnegative for all $n$.  In the process, we show how each of the resulting $\Pi$ can be obtained from a theorem  or conjecture involving more general sets of patterns.  In particular, we prove results concerning symmetries, shuffles, and Knuth classes, as well as pointing out a relationship with the arc permutations of Elizalde and Roichman.  Various conjectures and questions are mentioned throughout.
\end{abstract}

\section{Introduction}

Let $\fS_n$ be the symmetric group of all permutations of $[n]=\{1,\dots,n\}$.  We view the permutations in $\fS_n$ as sequences $\pi=\pi_1\dots\pi_n$.  Given any sequence $\si$ of $k$ distinct integers, its {\em standardization} is the permutation $\std\si\in\fS_k$ obtained by replacing the smallest element of $\si$ by $1$, the next smallest by $2$, and so forth.  We say that permutation $\si$ {\em contains permutation $\pi$ as a pattern} if there is a subsequence $\si'$ of (not necessarily consecutive) elements of $\si$ such that $\std\si'=\pi$.  For example, $\si=5132746$ contains $\pi=231$ because of the subsequence $\si'=574$.
Permutation $\si$ {\em avoids} $\pi$ if it does not contain $\pi$ as a pattern.  Given a set of permutations $\Pi$ we let 
$$
\fS_n(\Pi) =\{\si\in\fS_n \mid \text{$\si$ avoids every $\pi\in\Pi$}\}.
$$
For more information about permutation patterns, see the book of B\'ona~\cite{bon:cp}.  Our object is  to make a connection between the theory of patterns and the theory of quasisymmetric functions.  First we will  review some material concerning symmetric functions.  Details about symmetric functions and related combinatorics such as the Robinson-Schensted map can be found in the texts of Macdonald~\cite{mac:sfh}, Sagan~\cite{sag:sym}, or Stanley~\cite{sta:ec2}.

Let $\bx=\{x_1,x_2,\dots\}$ be a countably infinite set of commuting variables and consider the algebra of formal power series over the rationals $\bbQ[[\bx]]$.  Consider a monomial $m=x_{i_1}^{n_1}\dots x_{i_k}^{n_k}$.  The {\em degree} of $m$ is $\sum_i n_i$, and the {\em degree} of any $f\in\bbQ[[\bx]]$ is the maximum degree of a monomial in $f$ if the maximum exists or infinity otherwise.  We say that $f$ is {\em symmetric} if it is of bounded degree and invariant under permutation of variables.  
For example
$$
f = x_1^2 x_2 +  x_1 x_2^2 + x_1^2 x_3  + x_1 x_3^2 + x_2^2 x_3 + x_2 x_3^2 + \dots
$$
is symmetric.
Let $\Sym_n$ denote the vector space of symmetric functions homogeneous of degree $n$.  Bases for $\Sym_n$ are indexed by {\em partitions of $n$}, which are weakly decreasing sequences of positive integers $\la=(\la_1,\dots,\la_k)$ with $\sum_i \la_i=n$.
If $\la$ is a partition of $n$ then we write $\la\ptn n$ or $|\la|=n$.  The $\la_i$ are called the {\em parts} of $\la$.

We will be particularly interested in the important Schur basis for $\Sym_n$.  Recall that a partition $\la=(\la_1,\dots,\la_k)$ has an associated {\em Young diagram} consisting of $k$ left-justified rows of boxes with $\la_i$ boxes in row $i$.  We will write our diagrams in English notation with the first row on top and often make no distinction between a partition and its diagram.  Given $\la\ptn n$ then a {\em standard Young tableau (SYT)}, $P$, of shape $\la$ is obtained by filling the boxes bijectively with the elements of $[n]$ so that rows and columns increase.  In a {\em semistandard Young tableau (SSYT)}, $T$, of shape $\la$ the entries are positive integers distributed so that rows weakly increase and columns strictly increase.  
We write $\SYT(\la)$ or  $\SSYT(\la)$ for the set of standard or semistandard Young tableaux of shape $\la$, respectively. 
We also write $\sh P=\la$ or $\sh T=\la$ to indicate that $P$ or $T$ have shape $\la$.
The {\em Schur function} corresponding to $\la$ is
$$
s_\la=\sum_T \prod_{i\in T} x_i.
$$
For example, the set $\SSYT(2,1)$ consists of 
\beq
\label{21}
\begin{ytableau} 1&1\\ 2 \end{ytableau}\ , \hs{5pt}
\begin{ytableau} 1&2\\ 2 \end{ytableau}\ , \hs{5pt}
\begin{ytableau} 1&1\\ 3 \end{ytableau}\ , \hs{5pt}
\begin{ytableau} 1&3\\ 3 \end{ytableau}\ , \hs{5pt}
\dots,\quad
\begin{ytableau} 1&2\\ 3 \end{ytableau}\ , \hs{5pt}
\begin{ytableau} 1&3\\ 2 \end{ytableau}\ , \hs{5pt}
\begin{ytableau} 1&2\\ 4 \end{ytableau}\ , \hs{5pt}
\begin{ytableau} 1&4\\ 2 \end{ytableau}\ , \hs{5pt}
\dots
\eeq
resulting in 
$$
s_{(2,1)} = x_1^2 x_2 +  x_1 x_2^2 + x_1^2 x_3  + x_1 x_3^2 + \dots +  2 x_1 x_2 x_3 +  2 x_1 x_2 x_4 +\dots
$$
as the corresponding Schur function.
Call $f\in\Sym_n$ {\em Schur nonnegative} if its expansion in the $s_\la$ basis has nonnegative coefficients.

Quasisymmetric functions refine symmetric functions.  Call a formal power series $g$ of bounded degree {\em quasisymmetric} if any two monomials $x_{i_1}^{n_1}\dots x_{i_k}^{n_k}$ where $i_1<\dots<i_k$ and $x_{i_1}^{n_1}\dots x_{i_k}^{n_k}$ where $j_1<\dots<j_k$ have the same coefficient in $g$. By way of illustration,
$$
g=  x_1^2 x_2 + x_1^2 x_3   + x_1^2 x_4 + \dots + x_2^2 x_3 + x_2^2 x_4 +  \dots + x_3^2 x_4 + \dots
$$
is quasisymmetric, but not symmetric.  Clearly every symmetric function is quasisymmetric.  Quasisymmetric functions were introduced by Gessel~\cite{ges:mpi}  in his work on Stanley's theory of $P$-partitions.  They have since found many applications in both enumerative and algebraic combinatorics.  We denote by $\QSym_n$ the vector space of quasisymmetric functions homogeneous of degree $n$.  Bases for this vector space are indexed by {\em compositions of $n$}, which are sequences of positive integers  $\al=(\al_1,\dots,\al_k)$ summing to $n$.  We will use the notations $\al\comp n$ and $|\al|=n$ for compositions of $n$.  Greek letters near the beginning of the alphabet will be used for compositions while those near the middle will represent partitions.
There is also an important bijection between compositions of $n$ and subsets $S\sbe[n-1]$ given by
$$
(\al_1,\al_2,\al_3,\dots,\al_k)\mapsto \{\al_1,\ \al_1+\al_2,\ \al_1+\al_2+\al_3,\ \dots,\ \al_1+\al_2+\dots+\al_{k-1}\}.
$$
We will sometimes go back and forth between a composition and its associated set without mention.

The basis we will be using for $\QSym_n$ was considered in Gessel's original paper.  Given $S\sbe[n-1]$, the associated 
{\em fundamental quasisymmetric function} is
$$
F_S  = \sum x_{i_1} x_{i_2} \dots x_{i_n}
$$
where the sum is over all indices such that $i_1\le i_2\le \dots \le i_n$ and $i_s<i_{s+1}$ if $s\in S$.  For example, if $S=\{1\}\sbe[2]$ then the sum would be over all $x_i x_j x_k$ with $i<j\le k$, which gives
$$
F_{\{1\}} = x_1 x_2^2 + x_1 x_3^2 + x_2 x_3^2 + \dots + x_1 x_2 x_3 + x_1 x_2 x_4 + x_1 x_3 x_4 + x_2 x_3 x_4 + \dots
$$
in $\QSym_3$.   Since $s_\la$ is symmetric, and hence quasisymmetric, it can be expanded in the $F_S$ basis.  To do so, we need the  {\em descent set of an SYT}, $P$, which is
$$
\Des P =\{i \mid \text {$i+1$ is in a lower row than $i$ in $P$}\}.
$$
For example, if
$$
P=\raisebox{4mm}{\begin{ytableau} 1&2&5&9 \\ 3&4&7 \\ 6&8 \end{ytableau}}
$$
then $\Des P=\{2,5,7\}$.  Note that if $P$ has shape $\la\ptn n$ then $\Des P\sbe[n-1]$.
\bth[\cite{ges:mpi}]
\label{ges}
We have


\eqed{
s_\la = \sum_{P\in\SYT(\la)} F_{\Des P}.
}
\eth

Returning to $\la=(2,1)$, the two  elements of $\SYT(2,1)$  are the fifth and sixth displayed tableaux in~\ree{21} with descent sets $\{2\}$ and $\{1\}$ respectively.  Therefore $s_{(2.1)} = F_{\{1\}} + F_{\{2\}}$.

To combine permutation patterns and quasisymmetric functions, recall that the {\em descent set of a permutation}
$\pi=\pi_1\dots\pi_n$ is
$$
\Des\pi=\{i \mid \pi_i>\pi_{i+1}\}\sbe[n-1].
$$
For example, $\Des(35716824)=\{3,6\}$.
Now given a set of pattern permutations $\Pi$ we define the {\em pattern quasisymmetric function}
$$
Q_n(\Pi)=\sum_{\si\in\fS_n(\Pi)} F_{\Des\si}.
$$
The basic questions we wish to ask about these functions are
\ben
\item When is $Q_n(\Pi)$ symmetric for all $n$?
\item In that case, when is $Q_n(\Pi)$ Schur nonnegative for all $n$?
\een
Note that we are  asking  about symmetry or Schur nonnegativity not for a single function, but rather for an infinite family of functions.

\bfi
$$
\barr{l|l}
\Pi		& Q_n(\Pi) \text{ for $n\ge3$}\\
\hline
\emp 		&\dil \sum_\la  f^\la s_\la  \rule{0pt}{15pt}\\[20pt]
\{123\}	&\dil\sum_{\la_1<3}  f^\la s_\la\\[20pt]
\{321\}	&\dil\sum_{\la_1^t<3}  f^\la s_\la\\[20pt]
\{132,213\}; \{132,312\}; \{213,231\}; \{231,312\} &\dil\sum_{\text{$\la$ a hook}} s_\la\\[20pt]
\{123,132,312\}; \{123,213,231\}; \{123,231,312\} & s_{1^n} + s_{2,1^{n-2}}\\[10pt]
\{132,213,321\}; \{132,312,321\}; \{213,231,321\} & s_n + s_{n-1,1}\\[10pt]
\{132,213,231,312\}	& s_n + s_{1^n}\\[10pt]
\{123,132,213,231,312\} & s_{1^n} \\[10pt]
\{132,213,231,312,321\} & s_n 
\earr
$$
\capt{The $\Pi$ for Theorem~\ref{main} along with the Schur expansions of $Q_n(\Pi)$ where $\la\ptn n$
\label{S3}}
\efi

We start with $\Pi\sbe\fS_3$ and will prove the following theorem for which we will need some preliminaries.
If $\{123,321\}\sbe \Pi$ then, by the Erd\H{o}s-Szekeres Theorem~\cite{es:cpg}, $\fS_n(\Pi)=\emp$ for $n\ge5$, which explains the hypothesis on $\Pi$.  We use the notation $f^\la$  for the number of SYT of shape $\la$.  The {\em transpose} of $\la$ is the diagram $\la^t$ obtained by reflecting $\la$ in the main diagonal.  Then $\la_1^t$ is the number of boxes in the first column of $\la$, which is also written as $\ell(\la)$ and called the {\em length} of $\la$.
Also, a {\em hook} is a partition of the form $(a,1^b)$ for nonnegative integers $a,b$ where $1^b$ denotes the part $1$ repeated $b$ times.
When using a partition as a subscript we often omit the parentheses.  
\bth
\label{main}
Suppose $\{123,321\}\not\sbe \Pi\sbe\fS_3$.  The following are equivalent
\ben
\item   $Q_n(\Pi)$ is symmetric for all $n$.
\item   $Q_n(\Pi)$ is Schur nonnegative for all $n$.
\item   $\Pi$ is an entry in Table~\ref{S3}.
\een
\eth

We note that it is easy to compute examples to show that for any $\Pi$ not listed in Table~\ref{S3} we have $Q_n(\Pi)$ not being symmetric for some small value of $n$ depending on $\Pi$.  Therefore it suffices to show that the $\Pi$ in the table have the claimed Schur expansions.  In fact, we will show that these expansions are special cases of more general results or conjectures where $\Pi$ is not restricted to $\fS_3$.

The rest of this paper is structured as follows. In the next section we discuss what effect reversal and complementation have on $Q_n(\Pi)$ as well as using properties of the Robinson-Schensted correspondence to derive some of the results in Table~\ref{S3}.  In Section~\ref{shu} we show how the quasisymmetric function for a shuffle of two pattern sets can be computed in terms of the ones for each individual component.  Shuffles with increasing and decreasing permutations as well as full symmetric groups are used as examples.  We define partial shuffles in Section~\ref{psh} as certain shuffles where the increasing permutation has been removed.  We conjecture that in this case $Q_n(\Pi)$ has a nice Schur expansion and prove this in a special case.  Clearly if $\fS_n(\Pi)$ is a union of Knuth classes then $Q_n(\Pi)$ is symmetric and Schur nonnegative.  In Section~\ref{kc} we study when this can happen and, in particular, characterize the SYT such that the permutations  avoiding its Knuth class have this property.  We show in Section~\ref{ap} that avoiders of the arc permutations of Elizalde and Roichman~\cite{er:ap,er:sap} are in bijection with the permutations avoiding a certain set of shuffles.  We end with a section of comments and open questions.

\section{Symmetries and the Robinson-Schensted bijection}

If one views a permutation as a permutation matrix, then the dihedral group $D_4$ of the square acts on $\fS_n$.  We wish to investigate whether this action tells us anything about $Q_n(\Pi)$.  In particular, if this function is symmetric and one acts on $\Pi$ with a dihedral symmetry then is the new quasisymmetric function symmetric?  If so, what is its Schur expansion?  As usual, we apply a symmetry to a set by applying it to each member of the set.

First we note that if $Q_n(\Pi)$ is a symmetric function then $Q_n(\Pi^{-1})$ need not be.  For example, one can easily check that this is true if $n=3$ and $\Pi=\{132, 312\}$.  Therefore we can expect the dihedral symmetries preserving symmetry of $Q_n$ to form a subgroup of order at most $4$ in $D_4$.  We will now show that such a subgroup exists.

It is easy to describe two of the symmetries in $D_4$ directly in terms of the permutations in $\fS_n$.  The permutation $\pi=\pi_1\pi_2\dots\pi_n$ has {\em complement} $\pi^c = n+1-\pi_1, n+1-\pi_2,\dots,n+1-\pi_n$ where we have inserted commas between the elements for readability.  Its {\em reversal} is $\pi^r = \pi_n\pi_{n-1}\dots\pi_1$.  For example, if $\pi=35716824$ then $\pi^c = 64283175$ and $\pi^r = 42861753$.  Complementation and reversal are two of the reflections in $D_4$ and so $\pi^{rc}=\pi^{cr}$ is the permutation whose matrix is obtained by rotating the matrix of $\pi$ by $180^\circ$.  To state our first dihedral symmetry result, recall that  the transpose $\la^t$  of a Young diagram $\la$ is gotten by reflecting $\la$ in its main diagonal.  We use the same notation for tableaux.

\bpr
\label{c}
If $Q_n(\Pi)$ is symmetric then so is $Q_n(\Pi^c)$.  In particular, if $Q_n(\Pi)=\sum_\la c_\la s_\la$ for certain coefficients $c_\la$ then
$$
Q_n(\Pi^c)=\sum_\la c_\la s_{\la^t}.
$$
\epr
\bprf
Using Theorem~\ref{ges} we can write
$$
\sum_{\si\in\fS_n(\Pi)} F_{\Des\si} = Q_n(\Pi) = \sum_\la c_\la s_\la= \sum_\la c_\la \sum_{P\in\SYT(\la)} F_{\Des P}.
$$
Note that for any permutation $\Des \pi^c =[n-1]-\Des\pi$  and $\fS_n(\Pi^c)=(\fS_n(\Pi))^c$.  Also,
for any standard Young tableau $\Des P^t = [n-1]-\Des P$.  Using these facts and the previous displayed equation gives
$$
Q_n(\Pi^c)=\sum_{\si\in\fS_n(\Pi)} F_{[n-1]-\Des\si} = \sum_\la c_\la \sum_{P\in\SYT(\la)} F_{[n-1]-\Des P}=\sum_\la c_\la s_{\la^t}
$$
as desired.
\eprf

In order to deal with reversals, we will need some background about the Robinson-Schensted map~\cite{rob:rsg,sch:lid}.  This is a bijection
$$
\RS:\fS_n\ra\bigcup_{\la\ptn n} \SYT(\la)\times\SYT(\la).
$$
If $\RS(\pi)=(P,Q)$ then we will write $P=P(\pi)$ and $Q=Q(\pi)$ and call $P$ and $Q$ the {\em $P$-tableau} and {\em $Q$-tableau} of $\pi$, respectively.
If $\pi=\pi_1\dots\pi_n$ then $P$ is constructed using an operation called {\em insertion} that sequentially inserts $\pi_1,\dots,\pi_n$ to form $P$.  After the $k$th insertion, a $k$ is {\em placed} in $Q$ so that one maintains $\sh P=\sh Q$ at all times.  Although we are using $Q$ in both the notation $Q(\pi)$ and $Q_n(\Pi)$ the two different concepts should be distinguishable by context and the fact that the latter has a subscript while the former does not.  

We put an equivalence relation on $\fS_n$ by declaring that $\pi$ and $\si$ are {\em Knuth equivalent}, written $\pi\sim\si$, if $P(\pi)=P(\si)$.  Given an SYT $P$, the corresponding {\em Knuth class} is
$$
K(P) = \{\pi \mid P(\pi)=P\}.
$$
Given a partition $\la$ it will also be convenient to define an associated {\em Knuth aggregate} by
$$
K(\la) = \{\pi \mid \sh P(\pi) = \la\} =\bigcup_{\sh(P)=\la} K(P).
$$
We will need the following properties of the map $\RS$.
\bth
\label{RS}
Suppose $\RS(\pi)=(P,Q)$.
\ben
\item[(a)] $\Des\pi=\Des Q$.
\item[(b)] If $\sh P=\la$ then $\la_1$ is the length of a longest increasing subsequence of $\pi$.
\item[(c)] $P(\pi^r) = (P(\pi))^t$.
\item[(d)] $\RS(\pi^{-1}) = (Q,P)$. \hqed
\een
\eth

These results have interesting implications for our quasisymmetric functions.  To state them conveniently we will use the notation that for a permutation $\pi$ we let $F_\pi=F_{\Des\pi}$.  And for a set $\Pi$ of permutations we define $F_\Pi=\sum_{\pi\in\Pi} F_\pi$.
\ble
\label{knuth}
Suppose $P\in\SYT(\la)$.
\ben
\item[(a)] $F_{K(P)} = s_\la$.
\item[(b)] $F_{K(\la)} = f^\la s_\la$.
\een
\ele
\bprf
Using the fact that $\RS$ is a bijection, Theorem~\ref{ges}, and Theorem~\ref{RS} (a) we have
$$
F_{K(P)} = \sum_{\pi\in K(P)} F_{\Des\pi} = \sum_{Q\in\SYT(\la)} F_{\Des Q} = s_\la
$$
which is (a).   For (b), we use (a) to write
$$
F_{K(\la)} = \sum_{\sh(P)=\la} F_{K(P)} = \sum_{\sh(P)=\la} s_\la = f^\la s_\la
$$
as desired.
\eprf

We will now prove that Proposition~\ref{c} continues to hold if complement is replaced by reversal.
\bpr
\label{r}
If $Q_n(\Pi)$ is symmetric then so is $Q_n(\Pi^r)$.  In particular, if $Q_n(\Pi)=\sum_\la c_\la s_\la$ for certain coefficients $c_\la$ then
$$
Q_n(\Pi^r)=\sum_\la c_\la s_{\la^t}.
$$
\epr
\bprf
For each partition $\la$ pick a tableau $P_\la\in\SYT(\la)$.    Then using Lemma~\ref{knuth} (a) we have
$$
\sum_{\si\in\fS_n(\Pi)} F_{\Des\si} = Q_n(\Pi) = \sum_\la c_\la s_\la= \sum_\la c_\la F_{K(P_\la)}.
$$
Combining this with Theorem~\ref{RS} (c) yields
$$
Q_n(\Pi^r)=\sum_{\si\in\fS_n(\Pi)} F_{\Des\si^r} = \sum_\la c_\la F_{K((P_\la)^t)}=\sum_\la c_\la s_{\la^t}
$$
which is what we wished to prove.
\eprf

As an immediate corollary of Propositions~\ref{c} and~\ref{r} we get the following.
\bpr
\label{cr}
If $Q_n(\Pi)$ is symmetric then so is $Q_n(\Pi^{cr})$.  In particular,  $Q_n(\Pi^{cr})=Q_n(\Pi)$.\hqed
\epr

As an application of these ideas, we will now verify the expansions in the first three and last three rows of Table~\ref{S3}.  For the first row, we use Lemma~\ref{knuth} (b) and the fact that $\RS$ is a bijection to write
$$
Q_n(\emp)=\sum_{\si\in\fS_n} F_\si = \sum_{\la\ptn n} F_{K(\la)} = \sum_{\la\ptn n} f^\la s_\la.
$$
The next two rows are special cases of the following result. In it and the sequel we eliminate the braces when writing out $Q_n(\Pi)$ for a specific $\Pi$.  We will also use the notation
$$
\io_k = 12\dots k
$$
and
$$
\de_k = k \dots 21
$$
for the increasing and decreasing permutations of length $k$.
\bpr
\label{iode}
For any $k\ge 1$ we have
$$
Q_n(\io_k) = \sum_{\la_1<k} f^\la s_\la
$$
and
$$
Q_n(\de_k) = \sum_{\la_1^t<k} f^\la s_\la.
$$
\epr
\bprf
By either Proposition~\ref{c} or Proposition~\ref{r}, it suffices to prove the first statement.  Note that $\si\in\fS_n(\io_k)$ if and only if the longest increasing subsequence of $\si$ has length less than $k$.  Now using Theorem~\ref{RS} (b), Lemma~\ref{knuth} (b) and the bijectivity of $\RS$ we obtain
$$
Q_n(\io_k) =\sum_{\si\in\fS_n(\io_k)} F_\si = \sum_{\la_1<k} F_{K(\la)} =  \sum_{\la_1<k} f^\la s_\la
$$
which is the desired result.
\eprf

As far as  the last three rows of Table~\ref{S3}, the reader will find it easy to prove the following result so the proof is omitted.
\bpr
For $\Pi=\fS_k-\{\io_k,\de_k\}$ we have
$$
Q_n(\Pi)=\case{\dil\sum_{\la\ptn n} f^\la s_\la}{for $n<k$,}{\rule{0pt}{15pt} s_n+s_{1^n}}{for $n\ge k$.}
$$
For $\Pi=\fS_k-\{\de_k\}$ we have
$$
Q_n(\Pi)=\case{\dil\sum_{\la\ptn n} f^\la s_\la}{for $n<k$,}{\rule{0pt}{15pt} s_{1^n}}{for $n\ge k$.}
$$
For $\Pi=\fS_k-\{\io_k\}$ we have

\vs{10pt}
$\rule{1ex}{0ex}\hfill{  Q_n(\Pi)=\case{\dil\sum_{\la\ptn n} f^\la s_\la}{for $n<k$,}{\rule{0pt}{15pt} s_n}{for $n\ge k$.}  }\hfill\raisebox{-15pt}{\qed}$
\epr

\section{Shuffles}
\label{shu}

In this section we will show that several of the entries in Table~\ref{S3} can be explained using shuffles of permutations.  It turns out that under very general conditions, shuffling preserves $Q_n$ being symmetric, and perhaps Schur nonnegative as well.

If $\pi=\pi_1\dots\pi_m$ and $\si=\si_1\dots\si_n$ are sequences on distinct elements on disjoint alphabets then their {\em shuffle set} is
$$
\pi\shu\si =\{\tau=\tau_1\dots\tau_{m+n} \mid \text{$\pi$ and $\si$ are subsequences of $\tau$}\}.
$$
For example
$$
13\shu 42 = \{1342, 1432, 1423, 4132, 4123, 4213\}.
$$
If $\pi\in\fS_m$ and $\si\in\fS_n$ then we define $\pi\shu\si = \pi\shu(\si+m)$ where $\si+m$ is the sequence obtained by adding $m$ to each element of $\si$.  To illustrate
$$
12\shu 21 = 12 \shu 43 =\{1243, 1423, 1432, 4123, 4132, 4312\}.
$$
We shuffle sets of permutations as expected, namely
$$
\Pi\shu\Pi' = \bigcup_{\pi\in\Pi,\pi'\in\Pi'} \pi\shu\pi'.
$$

To prove the next result, it will be useful to have a notation for the permutations in $\fS_n$ that contain a pattern in $\Pi$, which will be
$$
\ofS_n(\Pi) = \fS_n - \fS_n(\Pi).
$$
We will also use $s_1$ as shorthand for the Schur function $s_{(1)}$.
\bth
\label{shu:thm}
 For any sets of permutations $\Pi$ and $\Pi'$ and any $n\ge0$,
$$
Q_n(\Pi \shuffle \Pi') = Q_n(\Pi') + \sum_{k=0}^{n-1} Q_k(\Pi) [s_1 Q_{n-k-1}(\Pi') - Q_{n-k}(\Pi')].
$$
\eth
\bprf
We first show that
\beq
\label{ofs}
\ofS_n(\Pi \shuffle \Pi') = \bigcup_{k=1}^{n-1} \ofS_k(\Pi) \shuffle \ofS_{n-k}(\Pi').
\eeq
If $\tau\in\ofS_n(\Pi \shuffle \Pi')$ then $\tau$ contains $\pi\shu\pi'$ for some $\pi\in\Pi$ and $\pi'\in\Pi'$.  Let $\tau^a$ and $\tau^b$ be the copies of $\pi$ and $\pi'$, respectively, in $\tau$. This implies $\max\tau^a<\min\tau^b$.  Then the restriction of $\tau$ to the elements of $[k]$ where $k=\max\tau^a$ is an element of $\ofS_k(\Pi)$, similarly restricting to $[n]-[k]$ gives an element of $\ofS_{n-k}(\Pi')+k$.  Thus $\tau\in\ofS_k(\Pi) \shuffle \ofS_{n-k}(\Pi')$.  The reverse inclusion is proven similarly.

The {\em Malvenuto-Reutenauer algebra}~\cite{mr:dqf}, $\MR$, is the set of formal  $\bbQ$-linear combinations of permutations with product given by shuffle.  The map $\Phi:\MR\ra\QSym$ given by $\pi\mapsto F_\pi$ is a homomorphism.  Define $A_k= \ofS_k(\Pi) \shuffle \ofS_{n-k}(\Pi')$.  Applying $\Phi$ to both sides of the summation identity implied by~\ree{ofs} and then using the Principle of Inclusion-Exclusion gives
\beq
\label{Phi(shu)}
\Phi(\ofS_n(\Pi \shuffle \Pi'))= \sum_{p=1}^{n-1} (-1)^{p-1} \sum_{1 \leq k_1 < \cdots < k_{p} \leq n-1} \Phi(A_{k_1} \cap \cdots \cap A_{k_p}).
\eeq
A proof similar to the one for~\ree{ofs} shows that if $k<\ell$ then 
\beq
\label{AcapA}
A_k \cap A_{\ell} = \ofS_k(\Pi) \shuffle \fS_{\ell-k} \shuffle \ofS_{n-\ell}(\Pi').
\eeq
It follows that if  $k_1 < \cdots < k_{p}$, then $A_{k_1} \cap \cdots \cap A_{k_p} = A_{k_1} \cap A_{k_{p}}$.   Applying this observation to~\ree{Phi(shu)} and then using the fact that $\sum_k (-1)^k \binom{n}{k}=\de_{n,0}$ gives
\begin{align*}
\Phi(\ofS_n(\Pi \shuffle \Pi')) 
&= \sum_{k=1}^{n-1} \Phi(A_k) + \sum_{1 \leq k < \ell \leq n-1} \Phi(A_k \cap A_\ell) \sum_{p=2}^{n-1} (-1)^{p-1} {\ell-k-1 \choose p-2}\\[10pt]
&= \sum_{k=1}^{n-1} \Phi(A_k) - \sum_{k=1}^{n-2} \Phi(A_k \cap A_{k+1}).
\end{align*}
Now using the definition of $A_k$, equation~\ree{AcapA}, and the fact that $\Phi$ is a homomorphism yields
\begin{align*}
\Phi(\ofS_n(\Pi \shuffle \Pi')) 
&= \sum_{k=1}^{n-1} \Phi(\ofS_k(\Pi) \shuffle \ofS_{n-k}(\Pi')) - \sum_{k=1}^{n-2} \Phi (\ofS_k(\Pi) \shuffle \fS_1 \shuffle \ofS_{n-k-1}(\Pi'))\\[10pt]
&= \sum_{k=1}^{n-1} \Phi(\ofS_k(\Pi)) \Phi(\ofS_{n-k}(\Pi')) - \sum_{k=1}^{n-2} \Phi (\ofS_k(\Pi))s_1\Phi(\ofS_{n-k-1}(\Pi')).
\end{align*}
Writing $\Phi(\ofS_n(\Pi))= \Phi(\fS_n) - Q_n(\Pi) = s_1^n - Q_n(\Pi)$, expanding the binomials, and rearranging the terms left after cancellation gives the theorem.
\eprf

We note that this theorem takes a nice form when expressed in terms of generating functions.  In particular, if one lets
$$
Q(\Pi) = \sum_{n=0}^{\infty} Q_n(\Pi)t^n
$$
then the previous result becomes the following.
\bco
For any sets of permutations $\Pi$ and $\Pi'$,

\eqed{
Q(\Pi \shuffle \Pi') = Q(\Pi) + Q(\Pi') + (ts_1 - 1)Q(\Pi)Q(\Pi').
}
\eco

We also have the following immediate corollary of Theorem~\ref{shu:thm}
\bco
For any sets of permutations $\Pi$ and $\Pi'$, if $Q_n(\Pi)$ and $Q_n(\Pi')$ are symmetric for all $n$ then then same is true of $Q_n(\Pi\shu\Pi')$.\hqed
\eco

Unfortunately, the theorem does not show that Schur nonnegativity is preserved.  However computer evidence supports this conjecture.
\bcon
\label{shu:non}
For any sets of permutations $\Pi$ and $\Pi'$, if $Q_n(\Pi)$ and $Q_n(\Pi')$ are Schur nonnegative for all $n$ then then same is true of $Q_n(\Pi\shu\Pi')$.
\econ

Although we can not do so in general, we can still derive Schur nonnegativity under certain circumstances.
\begin{lem}
\label{diff}
Suppose that, for all $n\ge0$
$$
G_n= \sum_{\la\ptn n} c_\la s_\la
$$
for certain constants $c_\la$.  Define, for all $n\ge1$, a symmetric function $G_n'$ and constants $d_\la$ by
$$
G_n'=s_1 G_{n-1} - G_n = \sum_{\la\ptn n} d_\la s_\la.
$$
\ben
\item[(a)]
We have
$$
d_\la = \left( \sum_{\la^-} c_{\la^-}  \right) - c_\la
$$
where $\la^-$ ranges over all diagrams obtained by removing a box from the diagram of $\la$
\item[(b)]
If $c_\la \le  \sum_{\la^-} c_{\la^-} $ for all $\la\ptn n$ then $G_n'$ is Schur-nonnegative.
\een
\end{lem}
\bprf
Clearly (ii) follows from (i).  And (i) itself follows easily from the Pieri rule for multiplying Schur functions.
\eprf

As an application, we will consider shuffles with $\io_k$, $\de_k$,  and $\fS_k$.
\begin{cor}
If $Q_n(\Pi)$ is Schur-nonnegative for all $n\ge0$ then so are the following: $Q_n(\Pi\shuffle\io_k)$, $Q_n(\Pi\shuffle\de_k)$, and $Q_n(\Pi\shuffle S_k)$.
\end{cor}
\bprf
We will only prove first statement as the others are similar. 
Recall that the Littlewood-Richardson rule expresses the product of two Schur functions as a non-negative linear combination of Schur functions.
Then in view of Theorem~\ref{shu:thm}, it suffices to show that condition (ii) of Lemma~\ref{diff} is satisfied where $G_n=Q_n(\io_k)$.  
By Proposition~\ref{iode}, if $\la_1\ge k$ then $c_\la=0$ and so the inequality is immediate.
On the other hand, if $\la_1<k$ then the same is true of all $\la^-$.  In this case
$$
c_\la=f^\la=\sum_{\la^-} f^{\la^-} =\sum_{\la^-} c_{\la^-}
$$
and we are done.
\eprf 

Now we can verify another four entries in Table~\ref{S3}.
\bth
For $n\ge2$ we have
$$
Q_n(123,132,312)=Q_n(123,213,231)= s_{1^n} + s_{2,1^{n-2}}
$$
and
$$
Q_n(132,312,321)=Q_n(213,231,321)=s_n + s_{n-1,1}.
$$
\eth
\bprf
By Propositions~\ref{c} and~\ref{r}, we only need to prove the statement for one of the four sets of permutations.  Consider 
$\{123,132,312\}=\{12\}\shu\{1\}$ and suppose $n\ge1$.  We have $Q_n(1)=\de_{n,1}$ where we are using the Kronecker delta.
Then for any $\Pi$,  Theorem~\ref{shu:thm} gives 
$Q_n(\Pi\shu 1) = s_1 Q_{n-1}(\Pi)$.  In particular, using Proposition~\ref{iode} and the Pieri rule
$$
Q_n(12\shu 1) = s_1 Q_{n-1}(12) = s_1 s_{1^{n-1}}=s_{1^n} + s_{2,1^{n-2}}
$$
as desired.
\eprf

\section{Partial shuffles}
\label{psh}

We now wish to study and generalize $Q_n(132,213)$ (and, by symmetry, $Q_n(231,312)$).  We will do this using a new concept that we call a partial shuffle. First, however, we will derive the quasisymmetric function for $\{132,213\}$ itself.

There is a well-known characterization of the permutations in $\fS_n(132,213)$; see, for example, the paper of Simion and Schmidt~\cite{ss:rp}.  A permutation in $\fS_n$ is {\em reverse layered} if it is of the form
\beq
\label{rlayer}
\pi = a, a+1,\dots, n, b, b+1, \dots, a-1, c, c+1, \dots, b-1, \dots
\eeq
For certain $a>b>c>\dots>0$.  We are using the term ``reverse" because in this case $\pi^r$ is what is usually called a {\em layered permutation}.  
Also, let
$$
\cH_n=\{\la\ptn n \mid \text{$\la$ is a hook}\}
$$
and $\SYT(\cH_n)$ be the set of tableaux $P$ so that $\sh(P) \in \cH_n$.

\bpr
\label{132,213}
We have
$$
Q_n(132,213)=Q_n(231,312)=\sum_{\la\in\cH_n} s_\la
$$
\epr
\bprf
By symmetry, it suffices to prove that $Q_n(132,213)$ is equal to the sum.
From the description~\ree{rlayer} it is clear that the map $\pi\mapsto\Des\pi$ gives a bijection $\fS_n(132,213)\ra 2^{[n-1]}$.  Thus
$$
Q_n(132,213)=\sum_{S\sbe[n-1]} F_S.
$$

We will be done by Lemma~\ref{knuth}  (a) if we can show that the map $P\mapsto\Des P$ is a bijection $\SYT(\cH_n)\ra 2^{[n-1]}$.   But is is easy to see that this map has an inverse.  In particular, given $S\sbe[n-1]$ we construct the hook tableau whose first column has elements $\{1\}\cup (S+1)$ where $S+1\sbs[n]$ is obtained by adding one to each element of $S$.
\eprf

To try and generalize the previous result, note that $\{132,213\}=(13\shu 2) - \{123\}$.  
We will put a hat on an element of a permutation to indicate the sequence obtained by removing that element.  For example, $3\hat{2}41=341$.
Now define the {\em partial shuffle}
$$
(12\dots\widehat{n-1}n)\pshu (n-1) = [(12\dots\widehat{n-1}n)\shu (n-1)] - \{\io_n\}.
$$
To illustrate $13\pshu 2 = \{132,213\}$ and $124\pshu 3 = \{1243, 1324, 3124\}$.  

For the generalization of hook diagrams, let $(i,j)$ denote the box of a Young diagram in row $i$ and column $j$.  
And if $P$ is a Young tableau then $P_{i,j}$ will be the entry of $P$ in box $(i,j)$.
We will consider the set of of {\em fattened hooks}
$$
\cH_{n,j} = \{\la\ptn n \mid (2,j)\not\in\la\}.
$$
Note $\cH_{n,2}=\cH_n$. The following conjecture has Proposition~\ref{132,213} as the  special case $j=3$.
\bcon
\label{pshu}
For $j\ge3$ we have
$$
Q_n(  (12\dots\widehat{j-1}j)\pshu (j-1) )=\sum_{\la\in\cH_{n,j-1}} f^{\lab} s_\la
$$
where  $\lab$ is $\la$ with $\la_1$ replaced by $\min\{\la_1,j-2\}$.
\econ

While we have not been able to prove this conjecture, we have been able to make some progress on the case $j=4$.
It follows from Theorem~\ref{RS} (b) and (c)  that if $\pi\in\Av(\de_m)$ then $P(\pi)$ has  less than $m$ rows.  In general it is {\em not} true that if $Q_n(\Pi)=\sum_\la c_\la s_\la$ then we have $Q_n(\Pi\cup \{\de_m\})$ is the same sum restricted to $\la$ with $\ell(\la)<m$.  As an example from Theorem~\ref{main}, 
$Q_n(231,312)$ is Schur non-negative but $Q_n(231,312,321)$ is not even symmetric for all $n$.  However, this property seems to be enjoyed in the context of partial shuffles.
\bcon
\label{pshum}
For $j\ge3$ and $m\ge 2$ we have
$$
Q_n(  (12\dots\widehat{j-1}j)\pshu (j-1) \cup\{\de_m\})=\sum_\la f^{\lab} s_\la
$$
where the sum is over all $\la\in\cH_{n,j-1}$ such that $\ell(\la)<m$, and $\lab$ is $\la$ with $\la_1$ replaced by $\min\{\la_1,j-2\}$.
\econ
Note that Conjecture~\ref{pshum} implies Conjecture~\ref{pshu} by letting $m\ra\infty$.  

We can give a proof of the  case $j=3$ of Conjecture~\ref{pshum}.  In fact, a slight modification of the proof of Proposition~\ref{132,213} will work since the bijections used there can be restricted to layered permutations with fewer than $m$ layers and hook tableaux with fewer than $m$ rows, both mapping onto the subsets of $[n-1]$ with fewer than $m$ elements.  This proves the following.
\bpr
\label{132,321,delta}
For $m\ge2$ we have
$$
Q_n(132,213,\de_m)=\sum_\la f^{\lab} s_\la
$$
where the sum is over all $\la\in\cH_n$ such that $\ell(\la)<m$, and $\lab$ is $\la$ with $\la_1$ replaced by $1$.  In particular
$$
Q_n(132,213,321)=s_n + s_{n-1,1}
$$
and by symmetry

\eqed{
 Q_n(123,231,312)=s_{1^n} + s_{2,1^{n-2}}.
}
\epr

We will now prove Conjecture~\ref{pshum} in the special case $j=m=4$, that is, we will consider 
\beq
\label{124pshu3,4321}
\Pi=\{1243, 1324, 3124, 4321\}.
\eeq
Given $\si\in\fS_n$ we will write its  descent set in increasing order $\Des\si=\{d_1<d_2<\dots<d_\ell\}$.   We also let $d_0=0$ and $d_{\ell+1}=n$ when convenient.  Let $\al(\si)$ denote the corresponding {\em descent composition} so that $\al_i = d_i-d_{i-1}$ for $1\le i \le \ell+1$.  The  {\em increasing run} of $\si$ corresponding to $\al_i$ is the factor of $\si$ consisting of the elements with indices between $d_{i-1}+1$ and $d_i$ inclusive.  For example, if $\si=561342$ then $\Des\si=\{2,5\}$, $\al(\si)=(2,3,1)$ and the increasing runs of $\si$ are $56$, $134$, and $2$.  We will need the following description of the increasing runs of a $\si$ avoiding $1243$ and $3124$.  It is analagous to the reverse layered description of $\si\in\Av_n(132,213)$, which is equivalent to  the increasing runs of $\si$   being consecutive integers.
\ble
\label{run}
If $\si$ avoids $1243$ and $3124$ then every increasing run of $\si$ of length at least two has the form
$$
a,b,b+1,b+2,\dots,b+k
$$
for some $a<b$ and some $k\ge0$.
\ele
\bprf
The result is trivial if the the length of the run is two, so assume it has length at least three.  Let $d$ be the maximum (and hence, last) element of the run.  Suppose, towards a contradiction, that not every integer in the interval $[b,d]$ appears in this run.  In particular, suppose that there is $c$ with $b<c<d$ that is not in the run.   Now either $c$ comes before the run so that $\si$ contains the subsequence $cabd$, or after the run so that $\si$ contains the subsequence $abdc$.  This
gives a contradiction in that $\si$ contains either a copy of $3124$ or $1243$.
\eprf

We now define two maps that in certain cases, because of the lemma, will be inverses.  The {\em contraction at index $j$} of $\si\in\fS_n$ is the permutation $c_j(\si)\in\fS_{n-1}$ obtained by removing $\si_j$ and then standardizing what remains.  The {\em expansion at index $j$} is the permutation $e_j(\si)\in\fS_{n+1}$ obtained by increasing all the elements of $\si$ greater than $\si_j$ by one and then inserting $\si_j+1$ directly after $\si_j$ in the result.  Continuing our example above with $\si= 561342$, we have $c_5(561342)=45132$ and also $e_4(45132)=561342=\si$.  To state the next result it will be convenient to have the notation.
$$
A(\al,n)=\{\si\in\fS_n(1243, 1324, 3124, 4321) \mid \al(\si)=\al\}.
$$
\bpr
\label{inverses}
Let $\al\comp n$ have a part $\al_i\ge3$ and let $j=\al_1+\al_2+\dots+\al_i$.  Also let $\al'\comp n-1$ be the composition obtained from $\al$ by replacing $\al_i$ by $\al_i-1$.  Then contraction and expansion restrict to maps $c_j:A(\al,n)\ra A(\al',n-1)$ and $e_{j-1}:A(\al',n-1)\ra A(\al,n)$ that are inverses of each other.
\epr
\bprf
We first show that $c_j$ is well defined.  Suppose $\si\in A(\al,n)$.  Since $\al_i\ge3$ we have, by Lemma~\ref{run}, that $\si_j=\si_{j-1}+1$.  It follows that removing $\si_j$ and standardizing will produce a descent at position $j-1$.  And since $c_j(\si)$ comes from standardizing a subpermutation of $\si$, it must still avoid the permutations in question.  Therefore $c_j(\si)\in A(\al',n-1)$ as desired.

The fact  that $e_{j-1}(\tau)=\si$ has the correct descent composition is similar to what was done for $c_j(\si)$.  To show that $\si$ still avoids the four patterns, assume to the contrary that it contains a copy $\ka$ of one of them.  Then $\ka$ must contain both $\si_{j-1}$ and $\si_j=\si_{j-1}+1$ because, if not,  then $\ka$ is also in $\tau$ (possibly obtained by replacing $\si_j$ by $\tau_{j-1}$).  Since $\si_j=\si_{j-1}+1$ these two elements of $\ka$ must correspond to consecutive,  increasing elements in the pattern, so the only option is that they become the one and two in either $1243$ or $3124$.  But $\al_i\ge3$ so that $\si_{j-1}$ is not the first element of its run.  Then replacing $\si_{j-1}\si_j$ by $\tau_{j-2}\tau_{j-1}$ gives a copy of the forbidden pattern in $\tau$, a contradiction.

We now show that the two maps are inverses.   For any permutation $\tau$ we have $c_j e_{j-1}(\tau)=\tau$.  On the other hand, if  $\si\in  A(\al,n)$ then, again using the fact that $\si_j=\si_{j-1}+1$, we must have $e_{j-1} c_j(\si)=\si$ because the element deleted by $c_j$ equals the element inserted by $e_{j-1}$.
\eprf

As for partitions, the number of parts of a composition $\al$ will be written $\ell(\al)$ and called the {\em length} of $\al$.  

\bco
\label{ell}
If $\si\in\fS_n(1243, 1324, 3124, 4321)$ then $\ell(\al(\si))\le 5$.
\eco
\bprf
Suppose, to the contrary, that $\ell(\al(\si))\ge 6$.  Removing all but the first $6$ increasing runs of $\si$ and standardizing we can assume $\ell(\al(\si))= 6$.  Now applying the bijections $c_j$ from the previous proposition for various values of $j$, we obtain a $\tau\in\fS_m(1243, 1324, 3124, 4321)$ that has no increasing run of length greater than $2$.  It follows that $m\le 12$.  But we have verified by computer that no such $\tau$ exists.
\eprf

We now have all the tools in place to compute $Q_n(\Pi)$ where $\Pi$ is given by~\ree{124pshu3,4321} and just need one more definition.   For compositions we will use the partial order $\al\le\be$ if $\ell(\al)=\ell(\be)$ and $\al_i\le\be_i$ for all $i$.  Note that we consider compositions of different lengths to be incomparable.

\bth
Let $\Pi=\{1243, 1324, 3124, 4321\}$.  Then
\beq
\label{lab}
Q_n(\Pi)=\sum_\la f^{\lab} s_\la
\eeq
where the sum is over all $\la\in\cH_{n,3}$ with $\ell(\la)\le 3$,  and $\lab$ is $\la$ with $\la_1$ replaced by $\min\{\la_1,2\}$.  In particular, for $n\ge6$ we have
\beq
\label{Qs}
Q_n(\Pi)=s_n + 2 s_{n-1,1} + 2 s_{n-2,2} + 3 s_{n-2,1,1} + 5 s_{n-3,2,1} + 5 s_{n-4,2,2}.
\eeq
\eth
\bprf
The ``in particular" follows from equation~\ree{lab} and the fact that when $n\ge6$ we always have $\la_1$ replaced by $2$ in $\lab$.

We will prove~\ree{lab} by induction on $n$, where we have verified it by computer for $n\le 10$.  Consider any nonempty set $A(\al,n)$ for $n\ge 11$.  By Corollary~\ref{ell}, the fact that this set is nonempty forces there to be a part $\al_i\ge 3$.  Applying the contraction bijections from Proposition~\ref{inverses}, we see that there is a $\be\comp 10$ such that $\#A(\be,10)=\#A(\al,n)$.

Now consider the sum on the right-hand side of~\ree{lab}, which, since $n\ge11$, reduces to~\ree{Qs}.  Using Theorem~\ref{ges}, we see that we want to show that $Q_n(\Pi)$ equals
\begin{align*}
&F_n + 2 \left[\sum_{\al\ge(1,1)} F_\al + \sum_{\al\ge(2,2)} F_\al  + \sum_{\al\ge(1,2,1)} F_\al\right ]
+  3 \sum_{\al\ge(1,1,1)} F_\al
\\[10pt]
&\quad + 5\left[  \sum_{\al\ge(2,2,1)} F_\al +  \sum_{\al\ge(2,2)} F_{1\al} +  \sum_{\al\ge(2,1,2)} F_\al
+ \sum_{\al\ge(1,2,1,1)} F_\al +  \sum_{\al\ge(1,1,2,1)} F_\al +  \sum_{\al\ge(2,2,2)} F_\al\right.
\\[10pt]
&\quad\left.  + \sum_{\al\ge(1,2,1,2)} F_\al +  \sum_{\al\ge(2,2,1)} F_{1\al} +  \sum_{\al\ge(2,1,2,1)} F_\al
+  \sum_{\al\ge(2,2,2,1)} F_\al +  \sum_{\al\ge(1,2,1,2,1)} F_\al +  \sum_{\al\ge(1,1,2,1,1)} F_\al
\right]
\end{align*}
where  $1\al$ is the composition obtained by concatenating $(1)$ and $\al$, and all subscripts of the $F$'s are compositions of $n$.  To finish the induction, it suffices to show that the   sums of fundamentals to which $\al$ contributes are the same as the ones to which $\be$ contributes where $\al$ and $\be$ are as in the previous paragraph.  But all of the lower bounds in the sums are  compositions with only ones and twos.  And applying $c_j$ reduces a part of size at least three to a part of size at least two and leaves all other parts the same.  Therefore this map does respect the expansion into fundamentals above and we are done.
\eprf

Note that~\ree{Qs} shows that the Schur expansion for this set of patterns stabilizes as $n\ra\infty$.  We will have more to say about stability in the next and last sections.

\section{Knuth classes}
\label{kc}

The final entries from Table~\ref{S3} that need to be explained are $\{132, 312\}$ and $\{213,231\}$.  The reader will recognize these as Knuth classes $K(P)$ for $P$ equal to
\beq
\label{12:3}
\begin{ytableau} 1&2\\ 3 \end{ytableau}
\eeq
and
\beq
\label{13:2}
\begin{ytableau} 1&3\\ 2 \end{ytableau}\ ,
\eeq
respectively.  The main result of this section will be a characterization of the $P$ for which $\fS_n(K(P))$ is a union of Knuth classes for all $n$.  For such $P$ it is clear that $Q_n(K(P))$ will be Schur nonnegative.  We will first deal with the special case when $P$ is one of the two tableaux above.

We will concentrate on $\{132, 312\}$ since the other $\Pi$ is just its reversal.  It follows easily from the inductive description of the $\si\in\fS_n(132,312)$ in~\cite{ss:rp} that they are exactly the permutations such that for all $j\ge2$ we have 
$$
\si_j=\case{(\min_{i<j}\{\si_i\})-1,}{or}{(\max_{i<j}\{\si_i\})+1.}{}
$$
Equivalently, for all $j\ge1$ the set $\{\si_1,\si_2,\dots,\si_j\}$ is an interval of integers. 
In the representation-theoretic proof of~\cite[Proposition 7.3]{er:ap} (see Section 10.2), Elizalde and Roichman mention a result closely related to the following.

\bpr
\label{132,312}
For $n\ge1$ we have
$$
Q_n(132,312)=Q_n(213,231)=\sum_\la s_\la
$$
where the sum is over all hook diagrams $\la$ with $n$ boxes.
\epr
\bprf
By reversal, it suffices to prove that $Q_n(132,312)$ is equal to the sum.  As in the proof of Proposition~\ref{132,213}, it suffices to show that $\Des$ restricts to a bijection $\fS_n(132,312)\ra 2^{[n-1]}$.  But given $J\sbe[n-1]$ there is a unique $\si$ avoiding $\{132,312\}$ with that descent set, which can be constructed as follows.  Arrange the integers $1,2,\dots,\#J$ in a decreasing sequence in the positions of $J+1$ (the set obtained by adding one to every element of $J$) of $\si$.  Then arrange the integers  $\#J+1,\#J+2,\dots,n$ as an increasing sequence in the remaining positions.  It is easy to check that the permutation thus constructed has $\Des\si=J$.
\eprf

We now consider general Knuth classes $K(P)$.  We first need to recall some general facts about these sets. Suppose we have  positive integers $a<b<c$ and a permutation $\si$ having a factor (subsequence of consecutive elements) $acb$ or $cab$.  Then we can perform a {\em Knuth move} on $\si$ by exchanging one factor for the other.  Alternatively, a Knuth move  can exchange  factors of the form $bac$  and $bca$.  Two permutations are {\em Knuth equivalent} if one can be transformed into the other by a sequence of Knuth moves.
\bth
\label{Kequiv}
If $P(\si)=P$ then $K(P)$ is the set of permutations Knuth equivalent to $\si$.\hqed
\eth

Call $\Pi$ {\em pattern-Knuth closed} if $\fS_n(\Pi)$ is a union of Knuth classes for all $n$.  Equivalently, the complement $\ofS_n(\Pi)$ is a union of Knuth classes for all $n$.  The following lemma is easy to prove directly from the definitions so its demonstration is omitted.
\ble
\label{pkc:rc}
If $\Pi$ is pattern-Knuth closed then so are $\Pi^c, \Pi^r,$ and $\Pi^{rc}$.\hqed
\ele

There is a sense in which this property is stable.
\ble
\label{pkc:M}
The set $\Pi$ is pattern-Knuth closed if and only if $\fS_n(\Pi)$ is a union of Knuth classes for $n\le M+1$ where $M$ is the maximum length of a permutation in $\Pi$.
\ele
\bprf
The forward direction is obvious so we will prove the converse.  It will be more convenient to show that $\ofS_n(\Pi)$ is a union of Knuth classes for $n>M+1$.  Suppose $\si\in\ofS_n(\Pi)$ contains a copy $\ka$ of some $\pi\in\Pi$.  By Theorem~\ref{Kequiv}, it suffices to show that the result $\si'$ of performing a Knuth move on $\si$ will still be in $\ofS_n(\Pi)$.  We will do this for a factor $acb$  where $a<b<c$ as the other possible Knuth moves can be handled similarly.  There are three cases.

If $\ka$ does not contain both $a,c$ then $\ka$ is still a subsequence of $\si'$  and we are done.  If $\ka$ contains all three elements, then $\si'$ contains  $\ka'$,  which is formed by replacing $acb$ by $cab$ in $\ka$.  Since our hypothesis implies that $\Pi$ is a union of Knuth classes, $\ka'$ is a copy of an element of $\Pi$ and we are again finished.  The final case is when $a,c$ are in $\ka$  but $b$ is not.  Let $\pi'$ be the standardization of the subsequence of $\si$ containing $\ka$ and $b$.  Then $\pi'\in\ofS_n(\Pi)$ for some $n\le M+1$.  It follows from the hypothesis in this direction that the permutation $\pi''$ obtained by switching the elements corresponding to $a,c$ in $\pi'$ still contains an element of $\Pi$.  Since $\si'$ contains a copy of $\pi''$ this finishes the proof.
\eprf

We now show that being pattern-Knuth closed is preserved by shuffling.
\ble
\label{pkc:shu}
If $A\sbe\fS_m$ and $B\sbe\fS_n$ are unions of Knuth classes then so is $A\shu B$.
\ele
\bprf
We will show that if $\al\in A$, $\be\in B$ and $\si\in \al\shu\be$ contains the factor $acb$ where $a<b<c$, then replacing this factor by $cab$ results in $\si'\in A\shu B$.  The proof for the other Knuth moves is similar.  There are three cases.

If $a\le m< c$  then $\al$ and $\be+m$ are still subwords of $\si'$ and so we are done.  If $c\le m$ then $acb$ is a factor of $\al$.  But $A$ is a union of Knuth classes so replacing this factor by $cab$ gives $\al'\in A$.  Thus $\si'\in\al'\shu\be \sbe A\shu B$. The only other possibility is $a>m$ in which case a similar argument to the one just given with $A$ replaced by $B$ completes the proof.
\eprf

Combining this lemma and equation~\ree{ofs} we immediately see that Conjecture~\ref{shu:non} is true in the case when $\Pi$ and $\Pi'$ are pattern-Knuth closed. In fact, we only need to assume that one of $\Pi$ and $\Pi'$ is pattern-Knuth closed.
\bpr
If $\Pi$ and $\Pi'$ are pattern-Knuth closed then so is $\Pi\shu\Pi'$. If $Q_n(\Pi)$ is Schur nonnegative for all $n$ and $\Pi'$ is pattern-Knuth closed, then $Q_n(\Pi\shu\Pi')$ is Schur nonnegative for all $n$.\hqed
\epr

\bprf
The first sentence follows from Lemma~\ref{pkc:shu} and equation~\ree{ofs}. Under the hypotheses of the second sentence, Theorem~\ref{shu:thm} implies that $Q_n(\Pi \shu \Pi')$ will be Schur nonnegative for all $n$ if the same is true of $s_1 Q_{n-1}(\Pi') - Q_n(\Pi')$. Recalling the homomorphism $\Phi : \MR \to \QSym$ sending $w \mapsto F_w$ from the proof of Theorem~\ref{shu:thm}, we have
$$
s_1 Q_{n-1}(\Pi') - Q_n(\Pi') = \Phi\left(\sum_{w \in A} w\right)
$$
where $A = (\fS_{n-1}(\Pi') \shu \iota_1) \setminus \fS_{n}(\Pi')$. But $A$ is a union of Knuth classes by Lemma~\ref{pkc:shu} and the assumption that $\Pi'$ is pattern-Knuth closed.
\eprf

We need one more  result before we prove the main theorem of this section. For  $J\sbe [n-1]$ we let
$$
D_J =\{\pi\in\fS_n \mid \Des\pi= J\}.
$$
Suppose that $\RS(\pi)=(P,Q)$.  Then, by Theorem~\ref{RS} (a), we have $\pi\in D_J$ if and only if $\Des Q= J$.  As with other operations, we apply the inverse operator to a set of permutations by applying it to each individual permutation.  It now follows from what we have just said and Theorem~\ref{RS} (d) that $\pi\in D_J^{-1}$ if and only if $\Des P = J$.  Then $D_J^{-1}$ is a union of Knuth classes.   Also, it follows easily from the definitions that $\Des(\pi^{-1})$ is the set of all $i$ such that $i+1$ appears to the left of $i$ in $\pi$.  These two observations are important for the proof of the following result.
\ble
\label{DJ-1}
For any $J$, the set  $D_J^{-1}$  is pattern Knuth closed.\hqed
\ele
\bprf
The proof is very similar to that of Lemma~\ref{pkc:M}.  We only need to take some care with the last case for which we use the same notation as in that demonstration.  Suppose first that the  elements of $\pi$ corresponding to $a,c$ in $\ka$ are not consecutive in value. From the discussion before the lemma, switching them will result in $\ka'$ whose standardization has the same inverse descent set as $\ka$, namely $J$.  Thus $\si'\in\ofS_n(D_J^{-1})$.  Now suppose that $a$ and $c$ standardize respectively  to $i$ and $i+1$ in $\pi$.  It follows that the other elements of $\ka$ are all less than $a$ or greater than $c$.  Now let $\ka'$ be $\ka$ with $c$ replaced by $b$.  Since $a<b<c$ we have that $\ka'$ also standardizes to $\pi$, and $\ka'$ is a subsequence of $\si'$, finishing this subcase.
\eprf

To state our principal result, we need one last set of definitions.  The {\em row superstandard} Young tableau of shape $\la$ is the SYT obtained by filling the first row with $1,2\dots,\la_1$; the next row with $\la_1+1,\la_1+2,\dots,\la_1+\la_2$; and so on.  Note that the tableau in~\ree{12:3} is row superstandard.  
A {\em column superstandard} Young tableau is defined similarly except that one fills the columns from left to right; see~\ree{13:2}.  A tableau of either type is just called {\em superstandard}

\bth
\label{pkc:main}
The class $K(P)$ is pattern-Knuth closed if and only if $P$ is a superstandard tableau of hook shape.
\eth

We will prove this in a sequence of propositions which deal with the various cases involved.  Note that because of Theorem~\ref{RS} (c) and Lemma~\ref{pkc:rc} we have $K(P)$ is pattern-Knuth closed if and only if $K(P^t)$ is.  Therefore we only need to prove closure, or lack thereof, for one of $P$ or $P^t$.  We also assume throughout that the shape of $P$ is $\lambda\ptn n$.

\bpr
\label{p:hook_closed}
If $P$ is superstandard of hook shape then $K(P)$ is pattern-Knuth closed.
\epr
\bprf
Suppose $\la=(n-k, 1^k)$.   By the remarks just before this proposition it suffices to prove this when $P$ is column  superstandard.  It is easy to see that in this case $K(P)=D_J^{-1}$ where $J=[k]$.  The result then follows from Lemma~\ref{DJ-1}.
\eprf

It will be convenient in our proofs to work with permutations and SYT having elements that are rational numbers, not just integers.  To this end, given an integer $a$ we let $a^+= a + 1/2$ and $a^- = a - 1/2$.  It will be important when comparing such permutations and tableaux that we always standardize them to have entries $[n]$ for some $n$. 
We will also need the {\em column reading word} of an SYT $P$, which is the permutation $\rho(P)$ obtained by recording the elements in each column of $P$ read bottom to top and then concatenating the sequences for the columns left to right.  It is easy to see that the insertion tableau of $\rho(P)$ is $P$.

\bpr
If $P$ is of hook shape but not superstandard then $K(P)$ is not pattern-Knuth closed.
\epr
\bprf
Our strategy, as suggested by Lemma~\ref{pkc:M}, will be to find an element  $\si\in\ofS_{n+1}(K(P))$ to which a Knuth move can be applied creating a permutation $\si'$ having no subsequence $\ka$   with $\std P(\ka)=P$. By transposition if necessary, we can assume that $n$ is in the first column of $P$.  Since $P$ is not superstandard, there is some $a>1$ in its first column such that $a+1$ is in the first row of $P$.  Let $a$ be the largest such element.  Let $b>a+1$ be the next element in $P$'s first column so that we have the following situation
$$
P=
\raisebox{15mm}{\ytableausetup
{boxsize=1.7em}
\begin{ytableau}
1 & \cdots & \scriptstyle a{+}1 & \cdots\\
\vdots\\
a\\
b\\
\vdots\\
n
\end{ytableau}}
$$
with $b+1$ being in the first column if $b<n$.

Let $\si$ be the permutation obtained by placing $a^+$ just before $b \in \rho(P)$ so that
$$
\si = n, \dots, a^+,  b, a, \dots
$$
which by construction is an element of $\ofS_{n+1}(\Pi)$.  Now exchange $a$ and $b$ to obtain 
$$
\si'= n,\dots,a^+,a,b,\dots
$$
Knuth equivalent to $\si$.
Let $\ka$ be any subsequence of $\si'$ with $|\ka|=n$ and suppose $k$ is the element removed to form $\ka$.  If $k=b$ then the first column of $P(\ka)$ contains $a^+$, which standardizes to $a+1$ so that $\std P(\ka)\neq P$.  If $k\neq b$ then when $b$ is inserted  it will enter the tableau in box $(1,2)$ since either $a$ or $a^+$ will be in box $(1,1)$ at that point.  When $P(\ka)$ is standardized, $b$ will either stay the same or be replaced by $b+1$.  But both those elements need to be in the first column, so at some point $b$ will be bumped out of the first row.  But then it will enter the second row in box $(2,2)$ because at least two of $1,a,a^+$ are in $\ka$ so that the element in box $(2,1)$ when $b$ is bumped will be $a$ or smaller.  Then there will be a $(2,2)$ entry in $P(\ka)$ forcing its shape not to be a hook, which completes our proof.
\eprf

Note that by performing transpositions, it suffices to prove Theorem~\ref{pkc:main} for tableaux $P$  satisfying $P_{2,1}=2$.

\bpr
Let $P$ be a standard tableau of non-hook shape with $P_{2,1}=2$.  If there exists an $i$ with $P_{i,1}>P_{i-1,2}$ then $K(P)$ is not pattern-Knuth closed.
\epr
\bprf
We will use the strategy and notation of the proof of the previous proposition.  We can write $P$ in the following form
$$
\begin{tikzpicture}[scale=.5]
\node at (-1,-2) {$P=$};
\draw (0,0) to (6,0) to (6,-2) to (5,-2) to (5,-3) to (3,-3) to (3,-4) to (2,-4) to (2,0);
\draw (0,-1) to (2,-1);
\draw (1,0) to (1,-2);
\draw (0,0) to (0,-5) to (1,-5) to (1,-4) to (2,-4) to (2,-2) to (0,-2);
\node at (.5,-.5) {$1$};
\node at (.5,-1.5) {$2$};
\node at (1.5,-.5) {$c$};
\node at (1.5,-1.5) {$d$};
\node at (1,-3) {$R$};
\node at (3.5,-1.5) {$S$};
\end{tikzpicture}.
$$
It is easy to see that the insertion tableau of $\pi= \rho(R), 2, d,  1, c, \rho(S)$ is $P$.  We consider
$$
\si = \rho(R), 2,  d,  2^+, 1, c, \rho(S),
$$
which is Knuth equivalent to
$$
\si' =  \rho(R), d, 2, 2^+, 1, c, \rho(S).
$$
Remove an element $k$ from $\si'$ to form $\ka$.  There are two cases.

Suppose first that $k\not\in S_1 \cup \{2,2^+,c\}$ where $S_1$ is the first row of tableau $S$.  It follows that $2, 2^+, c, S_1$ is an increasing subsequence of $\ka$ of length $\la_1+1$ where $\sh P = \la$. Thus, by Theorem~\ref{RS} (b), $P(\ka)$ has first row of length longer than $\la_1$ and so is not of the correct shape.  If $k\in S_1 \cup \{2,2^+,c\}$ then $k\neq 1,d$.  Let $i$ be an index with $P_{i,1}>P_{i-1,2}$ and note that $i\ge3$ so that $P_{i-1,2}\ge d$.  
It follows that $P_{\ell(\la),1}, \dots, P_{i,1}, P_{i-1,2}, \dots, P_{3,2}, d, x, 1$ is a decreasing subsequence of $\ka$ where $x$ is either $2$ or $2^+$.  But this subsequence has length $\ell(\la)+1$.  If follows from Theorem~\ref{RS} (b) and (c) that $P(\ka)$ will have first column of length longer than $\ell(\la)+1$ and so, again, will not have the right shape.  This finishes the proof.
\eprf

We can extend the definition of $P$ by letting $P_{i,j}=\infty$ for $i,j\ge1$ and $(i,j)$ outside $\la=\sh P$.    In this case, the proof of the previous proposition still goes through in the case where $P$'s first two columns are of equal length since $P_{\ell(\la)+1,1}=\infty>P_{\ell(\la),2}$.  We can now finish the proof of Theorem~\ref{pkc:main} with the following proposition.
\bpr
Let $P$ be of non-hook shape.  If $P_{i+1,1}<P_{i,2}$ for all $1\le i\le t$ where $t$ is the length of $P$'s second column, then $K(P)$ is not pattern-Knuth closed.
\epr
\bprf
We can write
\begin{center}
\begin{tikzpicture}[scale = .6]
\draw (0,0) to (8,0) to (8,-1) to (6,-1) to (6,-2) to (5,-2);
\draw (2,-2) to (0,-2) to (0,0);
\draw (5,-2) to (5,-3) to (4,-3) to (4,-4) to (0,-4) to (0,-2);
\draw (0,-3) to (2,-3);
\draw (2,-2) to (2,-4);
\draw (1,-2) to (1,-6) to (0,-6) to (0,-4);
\draw (2,0) to (2,-2);
\node at  (-1.3,-2) {$P=$};
\node at (1,-1) {$R$};
\node at (4,-2) {$S$};
\node at (.5,-5) {$C$};
\node at (.5,-2.5) {$a$};
\node at (.5,-3.5) {$b_0$};
\node at (1.5,-2.5) {$c$};
\node at (1.5,-3.5) {$e$};
\end{tikzpicture}
\end{center}
where $C$ is a single column.  Note that by the discussion just before this proposition, $C$ must be nonempty.  Let 
$\rho(C) = f_p \dots f_1 d_l \dots d_1 b_k \dots b_1$ where
$$
a < b_0 < b_1 < \dots < b_k < c < d_1 < \dots < d_l < e < f_1 < \dots < f_p.
$$
Note that either $k>0$ or $l>0$ since $C$ must contain an element less than $e$ to satisfy the hypothesis of this proposition.
Our proof splits into two cases.

\bfi
$$
\ytableausetup
{boxsize=1.7em}
\begin{ytableau}
R_1 &  R_2 & e\\
a &  c \\
b_0 & d_1 \\
b_1\\
\vdots\\
b_k\\
c^-\\
d_2\\
\vdots\\
d_l\\
f_1\\
\vdots\\
f_p
\end{ytableau}
$$
\capt{The insertion tableau for $\si'$ in Case 1.
\label{P(si')}}
\efi

\noindent{\bf Case 1:} $l > 0$.
Define
$$
\pi = \rho(C),  b_0,  e,  a, c,  \rho(R), \rho(S).
$$
It easy to check that $P(\pi) = P$.
Also let
$$
\si= f_p,\dots,f_1,d_l,\dots,d_2,c^-,d_1,b_k,\dots,b_0,  e,  a, c,  \rho(R), \rho(S)
$$
which is obtained from $\pi$ by adding $c^-$.  Apply a Knuth move to obtain
$$
\si'= f_p,\dots,f_1,d_l,\dots,d_2,c^-,b_k,d_1,b_{k-1},\dots,b_0,  e,  a, c,  \rho(R), \rho(S)
$$
Remove any $h$ from $\si'$ to obtain $\ka$.  We will show $P(\ka)\neq P$.

Consider the insertion tableau of $\si$ up to and including $\rho(R)$, which will have the form given in Figure~\ref{P(si')} where $R_1$ and $R_2$ are the two columns of $R$.  If $h$ comes from the first column of this tableau then its removal will cause all the entries to shift up by one, making the first column too short.  The only way to compensate for this is if the insertion of $\rho(S)$ causes $e$ to bump into the first column in row $r$ just below $b_0$.   But then 
$$
P(\ka)_{r,1}=e>d_1=P(\ka)_{r-1,2}
$$
which is a contradiction.  On the other hand, if $h$ comes from the second column of  Figure~\ref{P(si')} then $e$ will have to end up in the second column to preserve column lengths.  And if $h$ comes from a column of $S$ then the first two columns in the figure must remain undisturbed.  In either case the first column will contain more entries in the interval $[b_0,c]$ than $P$ does, which finishes this case.

{\bf Case 2:} $l=0$ and $k>0$.
Now we consider
$$
\pi= f_p,\dots,f_1,b_k,e,b_{k-1},\dots,b_0,   a, c,  \rho(R), \rho(S).
$$
It is easily checked that since $k>0$ we have $P(\pi)=P$.    Add $c^+$ to obtain
$$
\si= f_p,\dots,f_1,c^+,b_k,e,b_{k-1},\dots,b_0,   a, c,  \rho(R), \rho(S)
$$
and perform a Knuth move to yield
$$
\si'=f_p,\dots,f_1,c^+,e,b_k,b_{k-1},\dots,b_0,   a, c,  \rho(R), \rho(S).
$$
Define $h$, $\ka$, and $R_1$ as before with the goal of showing $P(\ka)\neq P$.

Note that there is a decreasing subsequence of $\si'$ of length $\ell(\la)+1$ consisting of the $f$'s, either $c^+$ or $e$, the $b$'s, $a$, and $\rho(R_1)$.
Therefore $h\in\{a,b_0,\dots,b_k,f_1,\dots,f_p\}\cup R_1$.  If $h$ is one of the $f$'s then $P(\ka)$ will have more elements in the interval $[a,e]$ in its first column than $P$.  If $h$ is any of the other possibilities, then $P(\ka)$ will have more elements greater than $c$ in its first column than $P$.  Either way, we have our final contradiction.
\eprf

We conclude this section with a two questions.  
First of all, we know from Lemma~\ref{knuth} Theorem~\ref{pkc:main} that when $P$ is superstandard of hook shape then $Q_n(\Pi)$ is  Schur nonnegative.
\begin{question}
Find a combinatorial interpretation for the coefficients in the Schur expansion for $Q_n(K(P))$ when $P$ is superstandard of hook shape.
\end{question}
It is also natrual to ask about  generalizing  Theorem~\ref{pkc:main}  to pairs of  tableaux.
\begin{question}
\label{c:2_knuth}
Let $P,Q$ be standard Young tableaux.  Find a necessary and sufficient condition for $K(P) \cup K(Q)$ to be pattern-Knuth closed.
 \end{question}

We note that if $\Pi,\Pi'$ are both pattern-Knuth closed then so is $\Pi\cup\Pi'$ since in this case $\ofS_n(\Pi\cup\Pi')$ is the union of the Knuth classes of $\ofS_n(\Pi)$ and $\ofS_n(\Pi')$. This together with Theorem~\ref{pkc:main} shows that if $P,Q$ are superstandard hooks then  $K(P) \cup K(Q)$ is pattern-Knuth.

\section{Arc permutations}
\label{ap}

Elizalde and Roichman introduced and studied arc permutations in~\cite{er:ap,er:sap}.  They are closely related to two of the permutation classes we have been considering and so we will be able to produce a bijection between the arc permutations and permutations avoiding a certain shuffle.

A permutation $\si\in\fS_n$ is in the set of {\em arc permutations}, $\cA_n$, if each prefix $\si_1\si_2\dots\si_i$ is an interval in the integers modulo $n$.  Equivalently, $\cA_n=\fS_n(\Pi_\cA)$ where
$$
\Pi_\cA=\{1324,1342,2413,2431,3124,3142,4213,4231\}.
$$
Recall that the permutations avoiding $\{132,312\}$ were those where each prefix was an interval of integers and a similar statement is true for suffixes when avoiding $\{213,231\}$.  Define the   {\em pin shaped} permutations to be
$$
\cP_n=\fS_n(132,312)\cup\fS_n(213,231)\sbe \cA_n.
$$
(Elizalde and Roichman call these permutations {\em unimodal} because their rotations are unimodal sequences, but we prefer to reserve unimodal for its original meaning.)  We will also be using  the complement $\cZ_n=\cA_n-\cP_n$.

In order to define the symmetric functions  to describe $Q_n(\Pi_\cA)$ we will need a few more definitions.  Consider the set of non-trivial hooks
$$
\cHb_n = \cH_n-\{(n), (1^n)\}
$$
as well as the diagrams one obtains from the elements of $\cHb_n$ by adding the $(2,2)$ box
$$
\cT_n = \{\la\cup (2,2) \mid \la\in\cHb_{n-1}\}.
$$
Define the corresponding generating functions $\Hb_n=\sum_{\la\in\cHb_n} s_\la$ and $T_n=\sum_{\la\in\cT_n} s_\la$.  The following result is a consequence of~\cite[Theorem 7.7]{er:ap}.
\bth
\label{er}
For all $n\ge0$ we have 

\eqed{Q_n(\Pi_\cA) = T_n + 2\Hb_n + s_{n} + s_{1^n}}
\eth

There is a shuffle class with the same expansion.  Let
$$
\Pi_\cS = \{1\}\shu \{132,312\}= \{1243,1423,2143,4123,2413,4213,2431,4231\}.
$$
\bpr
For all $n\ge0$ we have
$$
Q_n(\Pi_\cS)=Q_n(\Pi_\cA).
$$
\epr
\bprf
Applying Theorem~\ref{shu:thm} with $\Pi=\{1\}$ and $\Pi'=\{132,312\}$, then Proposition~\ref{132,312}, and finally the Pieri formula gives
$$
Q_n(\Pi_S) = Q_n(\Pi') + s_1 Q_{n-1}(\Pi') - Q_n(\Pi') = T_n + 2\Hb_n + s_{n} + s_{1^n}.
$$
Comparing this with the previous result completes the proof.
\eprf

Note that by reversal, complement and Proposition~\ref{132,213}, the previous proposition is true if $\Pi_\cS$ is replaced by either $\{1\}\shu\Pi$ or $\Pi\shu\{1\}$ for any
$$
\Pi\in \{\; \{132,312\}, \{132,213\}, \{213,231\},\{231,312\}\;\}.
$$
However, there is no dihedral symmetry relating any of these shuffles to $\Pi_\cS$.

Elizalde and Roichman~\cite[Section 7.4]{er:ap} gave a bijective proof of Theorem~\ref{er}.  In particular, they proved the following result.
\bth
For all $n\ge0$, there is an explicit bijection 

\eqed{\phi:\cZ_n\ra \bigcup_{\la\in\cT_n} \SYT(\la).}
\eth

As a consequence, we have the following.
\bco
For all $n\ge0$, there is an explicit bijection
$$
\psi:\fS_n(\Pi_\cA)\ra\fS_n(\Pi_\cS).
$$
\eco
\bprf
Note that by their descriptions in terms of prefixes and suffixes, we have 
$$
\fS_n(132,312)\cap\fS_n(213,231)=\{\io_n,\de_n\}.
$$
For $\si\in\fS_n(\Pi_\cA)$ we define $\psi(\si)$ as follows, using $Q(\pi)$ as the Robinson-Schensted recording tableau.
\ben
\item[(a)] If $\si\in\fS_n(132,312)$ then let $\psi(\si)=\si$.
\item[(b)] If $\si\in\cZ_n$ then let $\psi(\si)=\tau$ where $\tau\in\fS_n(\Pi_\cS)$ is the unique permutation with $Q(\tau)=\phi(\si)$.
\item[(c)] If $\si\in \fS_n(213,231)-\{\io_n,\de_n\}$ then let $\psi(\si)=\tau$ where $\tau\in\fS_n(\Pi_\cS)-\fS_n(132,312)$ is the unique permutation with $Q(\tau)=Q(\si)$.
\een
Note that (b) is well defined by Lemma~\ref{pkc:shu} and the fact that $Q_n(\cZ_n)$ is Schur multiplicity free.
In a  similar manner we see that (c) is well defined.   The check that this is a bijection is now easily done.
\eprf

The construction of the map $\psi$ is hardly as illuminating  as one would hope.

\begin{question}
Is there a direct description of a bijection between $\fS_n(\Pi_\cA)$ and $\fS_n(\Pi_\cS)$ on the level of permutations?
\end{question}

\section{Comments and open questions}
\label{coq}

We end with some comments and some questions that we hope the reader will be interested in answering.

\vs{15pt}

{\bf Symmetry vs nonnegativity.}
It is possible to construct $\Pi$ such that $Q_n(\Pi)$ is symmetric, but not Schur nonnegative.  In particular, consider the following set of permutations where we have enclosed certain elements in parentheses for readability
$$
\barr{l}
X_n=\{\io_n\} \cup \{2314\dots n,\ 12\dots(n-3)n(n-2)(n-1)\}\\[5pt]
\qquad \cup\{2134\dots n,\ 1324\dots n,\ \dots,\ 12\dots (n-2)n(n-1)\}\\[5pt]
\qquad \cup \{32145\dots n,\ 14325\dots n,\ \dots,\ 12\dots (n-3)n(n-1)(n-2)\}.
\earr
$$
Let $\Pi = \fS_4-X_4$.  Then one can verify that $\fS_n(\Pi)=X_n$ and so
$$
Q_n(\Pi)= s_n+2s_{n-1,1} + s_{n-2,1,1}-s_{n-2,2}.
$$

\vs{15pt}

%

{\bf Stability.}  The following is a natural question given the results we have proved such as Lemma~\ref{pkc:M}.
\begin{question}
Suppose $\Pi$ is a nonempty and $M$ is the maximum length of a permutation in $\Pi$.  Is there an $N$, a function of $M$, such that $Q_n(\Pi)$ symmetric for $n<N$ implies that it continues to be symmetric for $n\ge N$?  What about the same question with ``symmetric" replaced by ``Schur nonnegative?"
\end{question}

It is worth pointing out that the converse of these questions is false.  In particular, let $\la=(3,1,1)$ and 
$$
P=\raisebox{4mm}{\begin{ytableau}
1&2&4\\
3\\
5
\end{ytableau}}\ .
$$
Then Bloom and Sagan~[personal communication] have shown that  if $\Pi=K(\la)-K(P)$ then $Q_n(\Pi)$ is Schur nonnegative for $n\ge7$ but not even symmetric for $n=6$.

\vs{15pt}

{\bf Knuth classes.}  It would be interesting to determine when the union of Knuth classes is pattern-Knuth closed, generalizing Theorem~\ref{pkc:main}.  Of course, if $K$ and $L$ are pattern-Knuth closed then so is $K\cup L$.  In the case that one gets a pattern-Knuth closed class $\Pi$, one would also like to characterize the coefficients in the Schur expansion of $Q_n(\Pi)$.  Bloom and Sagan~[personal communication] have done this in the case where $\Pi=K(P)$ for a superstandard hook tableau $P$.

\vs{15pt}

{\bf Representation theory.} Adin and Roichman~\cite{ar:mcd} have developed  a way to connect certain subsets of permutations with representations of $\fS_n$.
Given $\al=(\al_1,\dots,\al_k)\comp n$ and $\pi\in\fS_n$ the {\em $\al$-decomposition of $\pi$} is the factorization
$\pi = \pi^{(1)} \dots \pi^{(k)}$ where $\#\pi^{(i)}=\al_i$ for all $i$. The {\em $\al$-descent set} and {\em $\al$-descent number} of $\pi$ are
$$
\Des_\al \pi = \bigcup_i \Des \pi_i \qmq{and} \des_\al \pi = \#\Des_\al \pi,
$$
respectively.   An integer sequence $a_1\dots a_p$ is {\em comodal} ({\em complement unimodal}) if, for some $m$, we have
$a_1>\dots > a_m < \dots < a_p$.  (Note that Adin and Roichman call such a sequence ``unimodal" but that is at variance with standard practice.)
Say that $\pi\in\fS_n$ is {\em $\al$-comodal} if each $\pi^{(i)}$ in its $\al$-decomposition is comodal.  If $\Pi$ is a set of permutations then let $\Pi_\al$ denote the set of $\al$-comodal permutations in $\Pi$.  Finally, call $\Pi\sbe\fS_n$ {\em fine} if there is an $\fS_n$-character $\chi$ such that for all $\al\comp n$
$$
\chi(\al) = \sum_{\pi\in\Pi_\al} (-1)^{\des_\al \pi}
$$
where $\chi(\al)$ is the value of $\chi$ on the conjugacy class indexed by $\al$. 

Adin and Roichman give a number of conditions equivalent to a set of permutations $\Pi\sbe\fS_n$ being fine, one of which is that $Q_n(\Pi)$ is symmetric and Schur nonnegative. 
They also state that the only known  examples of fine permutations are
\ben
\item permutations of given Coxeter length,
\item unions of conjugacy classes,
\item unions of Knuth classes, and
\item arc permutations.
\een
It is clear that the present work expands on this list.

\vs{15pt}

{\bf Other bases.}
A couple of the generating functions we have calculated are Schur $P$-functions.  For example, if $\Pi=\{132,312\}$ then $Q_n(\Pi)=P_n$ by Proposition~\ref{132,312}, and
$Q(\Pi_\cA)=P_{n-1,1}$ by Theorem~\ref{er}.
\begin{question}
For which $\Pi$ is $Q_n(\Pi)$ Schur $P$-nonnegative for all $n$?
\end{question}

One could also change the set of quasisymmetric functions being used in $\QSym_n$.  Let $\Th_{\Pk}$ be the peak fundamental quasisymmetric function associated with a  subset $\Pk\sbe[2,n-1]$.  Stembridge~\cite{ste:epp} developed a theory of enriched $P$-partitions ($P$ here is a partially ordered set), which allowed him to prove an analogue of Theorem~\ref{ges}   where $s_\la$ is replaced by a Schur $P$-function and the $F_{\Des P}$ are replaced by $\Th_{\Pk(\si)}$ for certain permutations $\si$, where $\Pk(\si)$ is the peak set of $\si$.  Then it is natural to define for a set of  permutations $\Pi$
$$
R_n(\Pi)=\sum_{\si\in\fS_n(\Pi)} \Th_{\Pk(\si)}
$$
and ask the question
\begin{question}
For which $\Pi$ is $R_n(\Pi)$ symmetric?  In that case when is it Schur $P$-nonnegative?
\end{question}

\vs{15pt}

{\em Acknowledgement.}  We would like to thank the following people for helpful discussions: Alex Burstein,  Vic Reiner, Yuval Roichman, Jon Stembridge, and  Damir Yeliussizov.
Joel Lewis suggested several ideas that were critical in arriving at our results in Section~\ref{kc}.



\end{document}